\documentclass[11pt,reqno]{amsart}

\usepackage{amsmath,amsthm}
\usepackage{amssymb}
\usepackage{euscript}

\numberwithin{equation}{subsection}

\newcommand{\sqsp}{\renewcommand{\baselinestretch}{1.2}\tiny\normalsize}

\newcommand{\relphantom[1]}{\mathrel{\phantom{#1}}}
\newtheorem{thm}[subsection]{Theorem}

\newtheorem{cor}[subsection]{Corollary}

\theoremstyle{definition}

\newcommand{\cat}[1]{{\EuScript #1}}

\newcommand{\cL}{\cat{L}}

\newcommand{\biglpren}{\bigl (}  % medium big (
\newcommand{\bigrpren}{\bigr )}  % medium big )
\newcommand{\llbrack}{\lbrack \lbrack}  % [[
\newcommand{\rrbrack}{\rbrack \rbrack}  % ]]
\newcommand{\Omegabar}{\overline{\Omega}}
\newcommand{\Omegatilde}{\widetilde{\Omega}}

\newcommand{\Deltabar}{\overline{\Delta}}

\newcommand{\Fbar}{\overline{F}}

\newcommand{\sumprime}{\sideset{}{'}\sum}  % \sum^\prime
\DeclareMathOperator{\Id}{Id}
\DeclareMathOperator{\Hom}{Hom}
\DeclareMathOperator{\Ob}{Ob}

\begin{document}
\title{Deformations of coalgebra morphisms}
\author{Donald Yau}

\begin{abstract}
An algebraic deformation theory of coalgebra morphisms is constructed.
\end{abstract}

%\subjclass[2000]{}
\address{Department of Mathematics, The Ohio State University at Newark, 1179 University Drive, Newark, OH 43055, USA}
\email{dyau@math.ohio-state.edu}
%\keywords{}
%\date{\today}
\maketitle
\sqsp

%%%================%%%
%%%                %%%
%%%  Introduction  %%%
%%%                %%%
%%%================%%%

\section{Introduction}
\label{sec:intro}

Algebraic deformation theory, as first described by Gerstenhaber \cite{ger1}, studies pertubations of algebraic structures using cohomology and obstruction theory.  Gerstenhaber's work has been extended in various directions. For example, Balavoine \cite{bal} describes deformations of any algebra over a quadratic operad.  In the same direction, Hinich \cite{hinich} studies deformations of algebras over a differential graded operad.  It is a natural problem to try to extend deformation theory to morphisms and, more generally, diagrams.

Deformations of morphisms are much harder to describe than that of the algebraic objects themselves.  Specifically, the difficulty arises when one tries to show that certain obstruction classes (obstructions to integration) are cocycles in the deformation complex.  Even in the case of an associative algebra morphism, proving this relies on a powerful result called the Cohomology Comparison Theorem (CCT) \cite{gs1,gs2,gs3}.  This Theorem says that the cohomology of a morphism, or a certain diagram, of associative algebras is isomorphic to the Hochschild cohomology of an auxiliary associative algebra.  This allows one to bypass the obstruction class issue by reducing the problem to the case of a single associative algebra.

The purpose of this paper is to establish the dual picture: A deformation theory of \emph{coalgebra morphisms}.  Deformations of coalgebras (not their morphisms) have been described by Gerstenhaber and Schack \cite{gs4}.  Their theory is essentially the same as that of the original theory \cite{ger1}, with Hochschild coalgebra cohomology in place of Hochschild cohomology for associative algebras.  We will build upon their work, regarding it as the absolute case, ours being the relative case.

There are two aspects of deformations of coalgebra morphisms that are different from the associative case.  First, the CCT requires a rather involved argument and has only been established for associative algebras.  Since we do not have a coalgebra version of the CCT, we deal with the obstruction class issue mentioned above differently, using instead an elementary, computational approach.  Second, working with coalgebras and their morphism seems to be even more conceptual and transparent than in the associative case.  In fact, much of our coalgebra morphism deformation theory is element-free.  In other words, although the Hochschild coalgebra cochain modules are of the form $\Hom(M,A^{\otimes n})$, we never need to pick elements in $M$ in our arguments.

It should be noted that deformations of Lie algebra morphisms have been studied by Nijenhuis-Richardson \cite{nr} and Fr\'eiger \cite{freiger}.  Also, recent work of Borisov \cite{borisov} sheds new light into the structure of the deformation complex of a morphism of associative algebras, showing that it is an $\cL_\infty$-algebra (``strongly homotopy Lie algebra").

\subsection*{Organization}
The next section is a preliminary one, in which Hochschild coalgebra cohomology is recalled. Section \ref{sec:mor def} introduces deformations of a coalgebra morphism and identifies infinitesimals as $2$-cocycles in the deformation complex $C^*_c(f)$ of a morphism of coalgebras (Theorem \ref{thm:inf}).  It ends with the Rigidity Theorem, which states that a morphism is rigid, provided that $H^2(C^*_c(f))$ is trivial (Theorem \ref{thm:rigidity}).  In section \ref{sec:ext}, the obstructions to extending a $2$-cocycle in $C^*_c(f)$ to a deformation of $f$ are identified. They are shown to be $3$-cocycles in $C^*_c(f)$ (Theorem \ref{thm:ob}).  In particular, such extensions automatically exist if $H^3(C^*_c(f))$ is trivial (Corollary \ref{cor:ext}).  Most of the arguments in this paper are contained in that section.

%%%%%%%%%%%%%%%%%%%%%%%%%%%%%%%%%%%%%%%%%
%%%%%%%%%%%%%%%%%%%%%%%%%%%%%%%%%%%%%%%%%
\section{Hochschild coalgebra cohomology}
\label{sec:coalg coh}

Throughout the rest of this paper, $K$ will denote a fixed ground field.  The following discussion of coalgebra cohomology is exactly dual to that of Hochschild cohomology of associative algebras.  See, for example, \cite{gs4}.

%%%%%%%%%%%%%%%%%%%%%%%
\subsection{Bicomodule}
A $K$-\emph{coalgebra} is a pair $(A, \Delta)$ (or just $A$) in which $A$ is a vector space over $K$ and $\Delta \colon A \to A\otimes A$ is a linear map, called \emph{comultiplication}, that is coassociative, in the sense that $(\Id \otimes \Delta) \Delta = (\Delta \otimes \Id) \Delta$.  When the ground field is understood, we will omit the reference to $K$.  If $A$ and $B$ are coalgebras, a \emph{morphism} $f \colon A \to B$ is a linear map that is compatible with the comultiplications, in the sense that $\Delta_B \circ f = (f \otimes f) \circ \Delta_A$.

An $A$-\emph{bicomodule} is a vector space $M$ together with a left action map, $\psi_l \colon M \to A \otimes M$, and a right action map, $\psi_r \colon M \to M \otimes A$, that make the usual diagrams commute.  For example, if $f \colon A \to B$ is a coalgebra morphism, then (the underlying vector space of) $A$ becomes a $B$-bicomodule via the structure maps:
   \[
   \begin{split}
   (f \otimes \Id_A) \circ \Delta_A & \colon A \to B \otimes A, \\
   (\Id_A \otimes f) \circ \Delta_A & \colon A \to A \otimes B.
   \end{split}
   \]
In this case, we say that $A$ is a $B$-bicomodule \emph{via} $f$.  One can also consider $A$ as an $A$-bicomodule with the structure maps $\psi_l = \psi_r = \Delta_A$.

%%%%%%%%%%%%%%%%%%%%%%%%%%%%%%%%%
\subsection{Coalgebra cohomology}
The \emph{Hochschild coalgebra cohomology} of a coalgebra $A$ with coefficients in an $A$-bicomodule $M$ is defined as follows.  For $n \geq 1$, the module of $n$-cochains is defined to be
   \[
   C_c^n(M,A) := \Hom_K(M, A^{\otimes n}),
   \]
with differential
   \[
   \begin{split}
   \delta_c \sigma
   & = (\Id_A \otimes \sigma) \circ \psi_l  \,+\, \sum_{i=1}^n (-1)^i \left(\Id_{A^{\otimes(i-1)}} \otimes \Delta \otimes \Id_{A^{\otimes(n-i)}}\right) \circ \sigma \\
   & \relphantom{} + (-1)^{n+1} (\sigma \otimes \Id_A) \circ \psi_r
   \end{split}
   \]
for $\sigma \in C^n_c(M,A)$.  Set $C^0_c(M,A) \equiv 0$.  The cohomology of the cochain complex $(C^*_c(M,A), \delta_c)$ is denoted by $H^*_c(M,A)$.

%%%%%%%%%%%%%%%%%%%%%%%%%%%%%%%%%%%%%%%%%
%%%%%%%%%%%%%%%%%%%%%%%%%%%%%%%%%%%%%%%%%
\section{Coalgebra morphism deformations}
\label{sec:mor def}

Fix a coalgebra morphism $f \colon A \to B$ once and for all.  Consider $A$ as a $B$-bicomodule via $f$ wherever appropriate.

The purpose of this section is to introduce deformations of a coalgebra morphism and discuss infinitesimals and rigidity.  All the assertions in this section are proved by essentially the same arguments as the ones in the associative case \cite{ger1,gs1,gs2,gs3}.  Therefore, we can safely omit the proofs.

%%%%%%%%%%%%%%%%%%%%%%%%%%%%%%%%
\subsection{Deformation complex}
\label{subsec:def complex}

For $n \geq 1$, define the module of $n$-cochains to be
   \[
   C^n_c(f) := C^n_c(A,A) \times C^n_c(B,B) \times C^{n-1}_c(A,B)
   \]
and the differential
   \[
   d_c \colon C^n_c(f) \to C^{n+1}_c(f)
   \]
by
   \[
   d_c(\xi; \pi; \varphi) =
   \left(\delta_c \xi;\, \delta_c \pi;\, \pi \circ f - f^{\otimes n} \circ \xi - \delta_c \varphi\right).
   \]
It is straightforward to check that $d_c d_c = 0$ \cite[p.\ 155]{gs3}.  The cochain complex $(C^*_c(f), d_c)$ (or just $C^*_c(f)$) is called the \emph{deformation complex of $f$}.  Its cohomology is denoted by $H^*_c(f)$.

Note that the signs in front of the terms $\pi \circ f$ and $f^{\otimes n} \circ \xi$ are different from their counterparts in the associative case \cite[p.\ 155, line 4]{gs3}.  This change of signs is needed in order to correctly identify infinitesimals as $2$-cocycles in the deformation complex (Theorem \ref{thm:inf}).

%%%%%%%%%%%%%%%%%%%%%%%%%%
\subsection{Deformation}
\label{subsec:deformation}

First, recall from \cite{gs4} that a \emph{deformation} of a coalgebra $A$ is a power series $\Delta_t = \sum_{n=0}^\infty \Delta_n t^n$ in which each $\Delta_n \in C^2_c(A,A)$ with $\Delta_0 = \Delta$, such that $\Delta_t$ is coassociative: $(\Id \otimes \Delta_t) \Delta_t = (\Delta_t \otimes \Id) \Delta_t$.  In particular, $\Delta_t$ gives a $K \llbrack t \rrbrack$-coalgebra structure on the module of power series $A \llbrack t \rrbrack$ with coefficients in $A$ that restricts to the original coalgebra structure on $A$ when setting $t = 0$.

With this in mind, we define a \emph{deformation of $f$} to be a power series $\Omega_t = \sum_{n=0}^\infty \omega_n t^n$, with each $\omega_n = (\Delta_{A,n};\, \Delta_{B,n};\, f_n) \in C^2_c(f)$, satisfying the following three statements:
   \begin{enumerate}
   \item $\Delta_{A,t} = \sum_{n=0}^\infty \Delta_{A,n}t^n$ is a deformation of $A$.
   \item $\Delta_{B,t} = \sum_{n=0}^\infty \Delta_{B,n}t^n$ is a deformation of $B$.
   \item $F_t = \sum_{n=0}^\infty f_nt^n \colon (A \llbrack t \rrbrack, \Delta_{A,t}) \to (B\llbrack t \rrbrack, \Delta_{B,t})$ is a $K\llbrack t \rrbrack$-coalgebra morphism with $f_0 = f$.
   \end{enumerate}
A deformation $\Omega_t$ will also be denoted by the triple $(\Delta_{A,t};\, \Delta_{B,t};\, F_t)$.

A \emph{formal isomorphism of $f$} is a power series $\Phi_t = \sum_{n=0}^\infty \phi_n t^n$ with each $\phi_n = (\phi_{A,n};\, \phi_{B,n}) \in C^1_c(f)$ and $\phi_0 = (\Id_A;\, \Id_B)$.

Suppose that $\Omegabar_t = (\Deltabar_{A,t};\, \Deltabar_{B,t};\, \Fbar_t)$ is also a deformation of $f$.  Then $\Omega_t$ and $\Omegabar_t$ are said to be \emph{equivalent} if and only if there exists a formal isomorphism $\Phi_t$ such that
   \begin{enumerate}
   \item $\Deltabar_{A,t} = \left(\Phi_{A,t} \otimes \Phi_{A,t}\right) \circ \Delta_{A,t} \circ \Phi_{A,t}^{-1}$,
   \item $\Deltabar_{B,t} = \left(\Phi_{B,t} \otimes \Phi_{B,t}\right) \circ \Delta_{B,t} \circ \Phi_{B,t}^{-1}$, and
   \item $\Fbar_t = \Phi_{B,t} \circ F_t \circ \Phi_{A,t}^{-1}$,
   \end{enumerate}
where $\Phi_{*,t} = \sum_{n=0}^\infty \phi_{*,n}t^n$ for $* = A, B$.

%%%%%%%%%%%%%%%%%%%%%%%%%%
\subsection{Infinitesimal}
\label{subsec:inf}

The linear coefficient $\omega_1 = (\Delta_{A,1};\, \Delta_{B,1};\, f_1) \in C^2_c(f)$ of a deformation $\Omega_t$ of $f$ is called the \emph{infinitesimal} of $\Omega_t$.  This element is more than just a $2$-cochain.

\begin{thm}
\label{thm:inf}
Let $\Omega_t = \sum_{n=0}^\infty \omega_n t^n$ be a deformation of $f$.  Then $\omega_1$ is a $2$-cocycle in $C^2_c(f)$ whose cohomology class is determined by the equivalence class of $\Omega_t$.  Moreover, if $\Omega_i = 0$ for $1 \leq i \leq l$, then $\omega_{l+1}$ is a $2$-cocycle in $C^2_c(f)$.
\end{thm}

%%%%%%%%%%%%%%%%%%%%%
\subsection{Rigidity}
\label{subsec:rigidity}

The \emph{trivial deformation of $f$} is the deformation $\Omega_t = \omega_0 = (\Delta_A;\, \Delta_B;\, f)$.  The morphism $f$ is said to be \emph{rigid} if and only if every one of its deformations is equivalent to the trivial deformation. The following cohomological criterion for rigidity is standard.

\begin{thm}
\label{thm:rigidity}
If $H^2_c(f)$ is trivial, then $f$ is rigid.
\end{thm}

%%%%%%%%%%%%%%%%%%%%%%%%%%%%%%%%%%%%%%%%%%%%%%%%
%%%%%%%%%%%%%%%%%%%%%%%%%%%%%%%%%%%%%%%%%%%%%%%%
\section{Extending $2$-cocycles to deformations}
\label{sec:ext}

In view of Theorem \ref{thm:inf}, a natural question is: Given a $2$-cocycle $\omega$ in $C^*_c(f)$, is there a deformation of $f$ with $\omega$ as its infinitesimal?  The purpose of this section is to identify the obstructions for the existence of such a deformation.  Following \cite{ger1}, if such a deformation exists, then $\omega$ is said to be \emph{integrable}.

Fix a positive integer $N$.   By a \emph{deformation of $f$ of order $N$}, we mean a polynomial $\Omega_t = \sum_{n=0}^N \omega_n t^n$ with each $\omega_n \in C^2_c(f)$ and $\omega_0 = (\Delta_A;\, \Delta_B;\, f)$, satisfying the definition of a deformation of $f$ modulo $t^{N+1}$.  In other words, for $X \in \lbrace A, B\rbrace$, $\Delta_{X,t} = \sum_{n=0}^N \Delta_{X,n}t^n$ defines a $K\lbrack t \rbrack/(t^{N+1})$-coalgebra structure on $X \lbrack t \rbrack/(t^{N+1})$ and $F_t = \sum_{n=0}^N f_nt^n$ is a $K\lbrack t \rbrack/(t^{N+1})$-coalgebra morphism.

To answer the integrability question, it suffices to consider the obstruction to extending $\Omega_t$ to a deformation of $f$ of order $N + 1$.  So let $\omega_{N+1} = (\Delta_{A,N+1};\, \Delta_{B,N+1};\, f_{N+1}) \in C^2_c(f)$ be a $2$-cochain and set
   \begin{equation}
   \label{eq:omegatilde}
   \Omegatilde_t := \Omega_t + \omega_{N+1}t^{N+1}.
   \end{equation}
Is $\Omegatilde_t$ a deformation of $f$ of order $N + 1$?  Since $\Omegatilde_t \equiv \Omega_t \pmod{t^{N+1}}$, it suffices to consider the coefficients of $t^{N+1}$ in the definition of a deformation of $f$.

To this end, consider the following cochains ($X \in \lbrace A, B\rbrace$):
   \begin{subequations}
   \begin{align}
   \Ob_X & = \sum_{i=1}^N \biglpren\left(\Delta_{X,i} \otimes \Id_X\right) \circ \Delta_{X,N+1-i} \,-\,
   \left(\Id_X \otimes \Delta_{X,i}\right) \circ \Delta_{X,N+1-i}\bigrpren, \notag \\
   \Ob_F & = \left(\sumprime \left(f_j \otimes f_k\right) \circ \Delta_{A,i}\right) \,-\, \sum_{i=1}^N \Delta_{B,N+1-i}\circ f_i, \label{eq:ObF} \\
   \intertext{where}
   \sumprime & = \sum_{\substack{i + j + k \,=\, N + 1 \\ 0 \,\leq\, i, j, k \,\leq\, N}}
   = \sum_{\substack{i + j \,=\, N + 1 \\ i, j > 0 \\ k = 0}} + \sum_{\substack{i + k \,=\, N + 1 \\ i, k > 0 \\ j = 0}} + \sum_{\substack{j + k \,=\, N + 1 \\ j, k > 0 \\ i = 0}} + \sum_{\substack{i + j + k \,=\, N + 1 \\ i, j, k > 0}}. \label{eq:sumprime}
   \end{align}
   \end{subequations}
From now on, integer indexes appearing in a summation are assumed non-negative, unless otherwise specified.  Let $\Ob_\Omega \in C^3_c(f)$ be the element
   \begin{equation}
   \label{eq:obomega}
   \Ob_\Omega = (\Ob_A;\, \Ob_B;\, \Ob_F).
   \end{equation}
Note that
   \[
   \Ob_X = \sum_{i=1}^N \Delta_{X,i} ~\bar{\circ}~ \Delta_{X,N+1-i} \in C^3_c(X,X)
   \]
for $X \in \lbrace A, B\rbrace$, where $\bar{\circ}$ is the ``comp" in \cite[p.\ 63]{gs4}.  A standard deformation theory argument \cite[p.\ 60]{gs4} shows that $\Ob_X$ is a $3$-coboundary if and only if  $\Delta_{X,t}$ extends to a $K\lbrack t \rbrack/(t^{N+2})$-coalgebra structure on $X \lbrack t \rbrack/(t^{N+2})$.  In this case, any $2$-cochain whose coboundary is $\Ob_X$ gives an extension.  An analogous argument applied to our setting yields the following result.

\begin{thm}
\label{thm:ext}
The polynomial $\Omegatilde_t$ is a deformation of $f$ of order $N + 1$ if and only if $\Ob_\Omega = d_c\omega_{N+1}$.
\end{thm}

It is also known \cite[Theorem 3]{gs4} that $\Ob_X$ is a $3$-cocycle.  We extend this statement to $\Ob_\Omega$.

\begin{thm}
\label{thm:ob}
The element $\Ob_\Omega \in C^3_c(f)$ is a $3$-cocycle.
\end{thm}

Before giving the proof of Theorem \ref{thm:ob}, let us record the following immediate consequence of the previous two Theorems.

\begin{cor}
\label{cor:ext}
If $H^3_c(f)$ is trivial, then every $2$-cocycle in $C^*_c(f)$ is integrable.
\end{cor}

\begin{proof}[Proof of Theorem \ref{thm:ob}]
Since $\Ob_A$ and $\Ob_B$ are $3$-cocycles, to show that $\Ob_\Omega$ is a $3$-cocycle, it suffices to show that
   \begin{equation}
   \label{eq:needtoshow}
   \delta_c \Ob_F \,-\, \Ob_B \circ f \,+\, f^{\otimes 3} \circ \Ob_A \,=\, 0.
   \end{equation}
To do this, first note that for a $2$-cochain $\psi \in C^2_c(A,B)$,  $\delta_c \psi$ is given by
   \[
   \begin{split}
   \delta_c \psi
   & = (\Id_B \otimes \psi) \circ (f \otimes \Id_A) \circ \Delta_A \,-\, (\Delta_B \otimes \Id_B) \circ \psi \\
   & \relphantom{} \relphantom{} + (\Id_B \otimes \Delta_B) \circ \psi \,-\, (\psi \otimes \Id_B) \circ (\Id_A \otimes f) \circ \Delta_A \\
   & = (f \otimes \psi) \circ \Delta_A - (\Delta_B \otimes \Id_B) \circ \psi + (\Id_B \otimes \Delta_B) \circ \psi - (\psi \otimes f) \circ \Delta_A.
   \end{split}
   \]
Applying this to $\Ob_F$ \eqref{eq:ObF}, we have
   \begin{subequations}
   \allowdisplaybreaks
   \begin{align}
   \delta_c \Ob_F & = \notag \\
   & \left( \sumprime f \otimes \lbrack (f_j \otimes f_k) \circ \Delta_{A,i}\rbrack \right) \circ \Delta_A \label{A} \\
   & \relphantom{} - \sum_{i=1}^N \left(f \otimes (\Delta_{B,N+1-i} \circ f_i)\right) \circ \Delta_A \label{B} \\
   & \relphantom{} - \sumprime (\Delta_B \otimes \Id_B) \circ (f_j \otimes f_k) \circ \Delta_{A,i} \label{C} \\
   & \relphantom{} + \sum_{i=1}^N (\Delta_B \otimes \Id_B) \circ \Delta_{B,N+1-i} \circ f_i \label{D} \\
   & \relphantom{} + \sumprime (\Id_B \otimes \Delta_B) \circ (f_j \otimes f_k) \circ \Delta_{A,i} \label{E} \\
   & \relphantom{} - \sum_{i=1}^N (\Id_B \otimes \Delta_B) \circ \Delta_{B,N+1-i} \circ f_i \label{F} \\
   & \relphantom{} - \left(\sumprime \lbrack (f_j \otimes f_k) \circ \Delta_{A,i} \rbrack \otimes f\right) \circ \Delta_A \label{G} \\
   & \relphantom{} + \sum_{i=1}^N \left( (\Delta_{B,N+1-i} \circ f_i) \otimes f \right) \circ \Delta_A. \label{H}
   \end{align}
   \end{subequations}
The terms \eqref{C}, \eqref{D}, and \eqref{E} need to be expanded.

For \eqref{C}, first note that, since $F_t = \sum_{n=0}^N f_nt^n$ is a $K \lbrack t \rbrack/(t^{N+1})$-coalgebra morphism, we have that
   \begin{equation}
   \label{eq:star}
   \sum_{i+j+k\,=\,n} (f_j \otimes f_k) \circ \Delta_{A,i}
   \,=\, \sum_{i=0}^n \Delta_{B,i} \circ f_{n-i}
   \end{equation}
for $0 \leq n \leq N$.  In particular, it follows that
   \begin{equation}
   \label{eq:starj}
   \Delta_B \circ f_j \,=\, \sum_{\alpha + \beta + \gamma \,=\, j} (f_\beta \otimes f_\gamma) \circ \Delta_{A,\alpha} \,-\, \sum_{\substack{\lambda + \mu\,=\, j \\ 1 \leq \mu \leq j}} \Delta_{B,\mu} \circ f_\lambda.
   \end{equation}
Putting \eqref{eq:starj} into \eqref{C}, we have that
   \begin{align*}
   \eqref{C}
   & = - \sumprime (\Delta_B \otimes \Id_B) \circ (f_j \otimes f_k) \circ \Delta_{A,i} \\
   & = - \sumprime \lbrack (\Delta_B \circ f_j) \otimes f_k \rbrack \circ \Delta_{A,i} \\
   & = - \sumprime_{\alpha + \beta + \gamma \,=\, j} \left\lbrack \left( (f_\beta \otimes f_\gamma) \circ \Delta_{A,\alpha}\right) \otimes f_k\right\rbrack \circ \Delta_{A,i} \\
   & \relphantom{} + \sumprime_{\substack{\lambda+\mu\,=\,j \\ 1 \leq \mu \leq j}} \left\lbrack (\Delta_{B,\mu} \circ f_\lambda) \otimes f_k \right\rbrack \circ \Delta_{A,i}.
   \end{align*}
The last two summations are given by
   \begin{subequations}
   \allowdisplaybreaks
   \begin{align}
   \sumprime_{\alpha + \beta + \gamma \,=\, j}
   & = \sum_{\substack{i + \alpha + \beta + \gamma \,=\, N+1 \\ i,\, \alpha + \beta + \gamma > 0 \\ k = 0}}
   + \sum_{\substack{i + k \,=\, N + 1 \\ i,k > 0 \\ \alpha = \beta = \gamma = 0}}
   + \sum_{\substack{\alpha + \beta + \gamma + k \,=\, N + 1 \\ k,\, \alpha + \beta + \gamma > 0 \\ i = 0}}
   + \sum_{\substack{i + \alpha + \beta + \gamma + k \,=\, N + 1 \\ i,\, k,\, \alpha + \beta + \gamma > 0}} \label{eq:sumprime1} \\
   \intertext{and}
   \sumprime_{\substack{\lambda + \mu \,=\, j \\ 1 \leq \mu \leq j}}
   & = \sum_{\substack{i + \lambda + \mu \,=\, N + 1 \\ i, \mu > 0 \\ k = 0}}
   + \sum_{\substack{\lambda + \mu + k \,=\, N + 1 \\ \mu, k > 0 \\ i = 0}}
   + \sum_{\substack{i + \lambda + \mu + k \,=\, N + 1 \\ i, \mu, k > 0}}. \label{eq:sumprime2}
   \end{align}
   \end{subequations}
Therefore, \eqref{C} can be written as
   \begin{subequations}
   \allowdisplaybreaks
   \begin{align}
   - \sumprime & \lbrack (\Delta_B \circ f_j) \otimes f_k \rbrack \circ \Delta_{A,i} =  \notag \\
   & - \sum_{\substack{i+\alpha \,=\, N+1 \\ i,\alpha > 0}} \left\lbrack \left((f \otimes f) \circ \Delta_{A,\alpha}\right) \otimes f\right\rbrack \circ \Delta_{A,i} \label{C111} \\
   & - \sum_{\substack{i+\alpha + \beta+\gamma\,=\,N+1 \\ i, \, \beta + \gamma >0}} \left\lbrack \left((f_\beta \otimes f_\gamma) \circ \Delta_{A,\alpha}\right) \otimes f\right\rbrack \circ \Delta{A,i} \label{C112} \\
   & - \sum_{\substack{i+k\,=\,N+1 \\ i,k>0}} \left\lbrack\left((f \otimes f)\circ \Delta_A\right) \otimes f_k \right\rbrack \Delta_{A,i} \label{C12} \\
   & - \sum_{\substack{\alpha+\beta+\gamma+k\,=\,N+1 \\ k, \alpha+\beta+\gamma > 0}} \left\lbrack \left((f_\beta \otimes f_\gamma) \circ \Delta_{A,\alpha}\right) \otimes f_k\right\rbrack \circ \Delta_A \label{C13} \\
   & - \sum_{\substack{i+\alpha+\beta+\gamma+k \,=\, N+1 \\ i, k, \alpha+\beta+\gamma > 0}} \left\lbrack \left( (f_\beta \otimes f_\gamma) \circ \Delta_{A,\alpha}\right) \otimes f_k\right\rbrack \circ \Delta_{A,i} \label{C14} \\
   & + \sumprime_{\substack{\lambda+\mu\,=\,j \\ 1 \leq \mu \leq j}} \left\lbrack (\Delta_{B,\mu} \circ f_\lambda) \otimes f_k \right\rbrack \circ \Delta_{A,i}. \label{C2}
   \end{align}
   \end{subequations}

A similar argument can be applied to \eqref{D}, which yields
   \begin{subequations}
   \allowdisplaybreaks
   \begin{align}
   \eqref{D}
   & \,=\, \sum_{i=1}^N (\Delta_B \otimes \Id_B) \circ \Delta_{B,N+1-i} \circ f_i \notag \\
   & \,=\, \sum_{i=1}^N (\Id_B \otimes \Delta_B) \circ \Delta_{B,N+1-i} \circ f_i \notag \\
   & \relphantom{} + \sum_{\substack{i+j+k\,=\,N+1 \\ i,k > 0}} (\Id_B \otimes \Delta_{B,k}) \circ \Delta_{B,j} \circ f_i \notag \\
   & \relphantom{} - \sum_{\substack{i+j+k\,=\,N+1 \\ i,k > 0}} (\Delta_{B,k} \otimes \Id_B) \circ \Delta_{B,j} \circ f_i \notag \\
   & \,=\, \sum_{i=1}^N (\Id_B \otimes \Delta_B) \circ \Delta_{B,N+1-i} \circ f_i \label{D1} \\
   & \relphantom{} + \sum_{\substack{k+\gamma \,=\, N + 1 \\ k,\gamma > 0}} \lbrack f \otimes (\Delta_{B,k} \circ f_\gamma)\rbrack \circ \Delta_A \label{D211} \\
   & \relphantom{} + \sumprime_{\substack{\lambda+\mu\,=\,k \\ 1 \leq \mu \leq k}} \lbrack f_j \otimes (\Delta_{B,\mu}\circ f_\lambda)\rbrack \circ \Delta_{A,i} \label{D212} \\
   & \relphantom{} - \sum_{\substack{j+k\,=\,N+1 \\ j,k>0}} (\Id_B \otimes \Delta_{B,k}) \circ \Delta_{B,j} \circ f \label{D22} \\
   & \relphantom{} - \sum_{i=1}^N \lbrack (\Delta_{B,N+1-i} \circ f_i) \otimes f\rbrack \circ \Delta_A \label{D31} \\
   & \relphantom{} - \sumprime_{\substack{\lambda+\mu\,=\,j \\ 1 \leq \mu \leq j}} \lbrack (\Delta_{B,\mu} \circ f_\lambda) \otimes f_k\rbrack \circ \Delta_{A,i} \label{D32} \\
   & \relphantom{} + \sum_{\substack{j+k\,=\,N+1 \\ j,k>0}} (\Delta_{B,k} \otimes \Id_B) \circ \Delta_{B,j} \circ f. \label{D33}
   \end{align}
   \end{subequations}

A similar argument when applied to \eqref{E} gives
   \begin{subequations}
   \allowdisplaybreaks
   \begin{align}
   \eqref{E}
   & = \sumprime \left(f_j \otimes (\Delta_B \circ f_k)\right) \circ \Delta_{A,i} \notag \\
   & = \sum_{\substack{i+j\,=\,N+1 \\ i,j>0}} \left\lbrack f_j \otimes \left((f\otimes f)\circ \Delta_A\right)\right\rbrack \circ \Delta_{A,i} \label{E11} \\
   & \relphantom{} + \sum_{\substack{i+\alpha\,=\,N+1 \\ i,\alpha>0}} \left\lbrack f \otimes \left((f \otimes f) \circ \Delta_{A,\alpha}\right)\right\rbrack \circ \Delta_{A,i} \label{E121} \\
   & \relphantom{} + \sum_{\substack{i+\alpha+\beta+\gamma\,=\,N+1 \\ i, \beta+\gamma > 0}} \left\lbrack f \otimes \left((f_\beta \otimes f_\gamma) \circ \Delta_{A,\alpha}\right)\right\rbrack \circ \Delta_{A,i} \label{E122} \\
   & \relphantom{} + \sum_{\substack{j+\alpha+\beta+\gamma\,=\,N+1 \\ j, \alpha+\beta+\gamma > 0}} \left\lbrack f_j \otimes \left((f_\beta \otimes f_\gamma) \circ \Delta_{A,\alpha}\right)\right\rbrack \circ \Delta_A \label{E13} \\
   & \relphantom{} + \sum_{\substack{i+j+\alpha+\beta+\gamma\,=\,N+1 \\ i, j, \alpha+\beta+\gamma > 0}} \left\lbrack f_j \otimes \left((f_\beta \otimes f_\gamma) \circ \Delta_{A,\alpha}\right)\right\rbrack \circ \Delta_{A,i} \label{E14} \\
   & \relphantom{} - \sumprime_{\substack{\lambda+\mu\,=\,k \\ 1 \leq \mu \leq k}} \left\lbrack f_j \otimes (\Delta_{B,\mu} \circ f_\lambda)\right\rbrack \circ \Delta_{A,i} \label{E2},
   \end{align}
   \end{subequations}
where the last sum is interpreted as in \eqref{eq:sumprime2} with the roles of $j$ and $k$ interchanged.

Now observe that each of the following sums is equal to $0$:
$-\Ob_B \circ f +$ \eqref{D22} $+$ \eqref{D33},
$f^{\otimes 3} \circ \Ob_A +$ \eqref{C111} $+$ \eqref{E121},
\eqref{D1} $+$ \eqref{F},
\eqref{D211} $+$ \eqref{B},
\eqref{D212} $+$ \eqref{E2},
\eqref{D31} $+$ \eqref{H},
\eqref{D32} $+$ \eqref{C2}.  It follows that
  \[
  \begin{split}
  \delta_c \Ob_F & \,-\, \Ob_B \circ f \,+\, f^{\otimes 3} \circ \Ob_A \\
  & = \sum \lbrace \lbrack f_\alpha \otimes ((f_\beta \otimes f_\gamma) \circ \Delta_{A,\mu})\rbrack \circ \Delta_{A,\lambda} \\
  & \relphantom{} \relphantom{} \relphantom{} - \lbrack ((f_\alpha \otimes f_\beta) \circ \Delta_{A,\mu}) \otimes f_\gamma\rbrack \circ \Delta_{A,\lambda}\rbrace.
  \end{split}
  \]
The sum on the right-hand side is taken over all $\alpha, \beta, \gamma, \lambda, \mu \geq 0$ such that:
   \begin{enumerate}
   \item $\alpha + \beta + \gamma + \lambda + \mu = N + 1$ with $1 \leq \alpha + \beta + \gamma \leq N$, or
   \item $\alpha + \beta = N + 1$ with $\alpha, \beta > 0$ and $\gamma = \lambda = \mu = 0$, or
   \item $\alpha + \gamma = N + 1$ with $\alpha, \gamma > 0$ and $\beta = \lambda = \mu = 0$, or
   \item $\beta + \gamma = N + 1$ with $\beta, \gamma > 0$ and $\alpha = \lambda = \mu = 0$, or
   \item $\alpha + \beta + \gamma = N + 1$ with $\alpha, \beta, \gamma > 0$ and $\lambda = \mu = 0$.
   \end{enumerate}
This sum is equal to $0$, since $\Delta_{A,t} = \sum_{i=0}^N \Delta_{A,i}t^i$ gives a $K \lbrack t \rbrack/(t^{N+1})$-coalgebra structure on $A \lbrack t \rbrack/(t^{N+1})$.

This finishes the proof of Theorem \ref{thm:ob}
\end{proof}

\section*{Acknowledgment}

The author thanks the referee for reading an earlier version of this paper and for pointing out the reference \cite{hinich}.

%%%==============%%%
%%%              %%%
%%%  References  %%%
%%%              %%%
%%%==============%%%

\end{document}